\documentclass[12pt]{amsart}

\usepackage{amsmath}
\usepackage{latexsym}
\usepackage{amssymb}
\usepackage[dvips]{graphicx}

\def\d{\bmath \delta}

\def\e{\bmath \epsilon}
\def\ee{\mathbf e}

\def\l{\bmath \lambda}

\def\mob{\mathfrak{m\ddot{o}b}}
\def\o{\omega}

\def\m{\mathbf m}
\def\t{{}^t\hskip-1pt}

\def\span{\operatorname{span}}

\def\x{\mathbf x}
\def\y{\mathbf y}

\def\n{\mathbf n}

\def\GL{\mathbf{GL}}

\def\0{\mathbf 0}
\def\B{\mathbf B}

\def\E{\mathbf E}

\def\g{\mathfrak g}
\def\gl{\mathfrak{gl}}
\def\h{\mathfrak h}
\def\H{\mathbf H}

\def\M{\mathcal M}
\def\Mob{\mathbf{M\ddot{o}b}}
\def\mob{\mathfrak{m\ddot{o}b}}
\def\Os{\mathbf S^{3,1}}
\def\R{\mathbf R} 

\def\SO{\mathbf{SO}}
\def\s2{\mathbf S^2}
\def\sb{\mathbf S}
\def\sph{\mathbf S^3}

\def\Tb{\mathbf T}
\def\bmath#1{\mbox{\boldmath$#1$}}

\def\Y{\mathbf Y}
\newtheorem*{congruence}{Congruence Theorem}
\newtheorem*{existence}{Existence Theorem}
\newtheorem{theorem}{Theorem}
\newtheorem{lemma}[theorem]{Lemma}

\theoremstyle{definition}
\newtheorem{definition}[theorem]{Definition}
\newtheorem{example}[theorem]{Example}

\theoremstyle{remark}

\begin{document}

\title{Dupin hypersurfaces in Lie sphere geometry}
\author{Gary R.\ Jensen}
\address{Department of Mathematics, Washington University in St. Louis}
\maketitle

\section*{Introduction}
The method of moving frames in Lie sphere geometry has produced significant
results in the classification of Dupin hypersurfaces in spheres, as
seen, for example in the papers by Pinkall
\cite{Pinkall}, Cecil and Chern \cite{CecChern}, Cecil and Jensen
\cite{CJ, CJ2}, Cecil, Chi, and Jensen \cite{DHwFPCII}.
What is the secret of its effectiveness?
The answer emerges in the classification of nonumbilic isoparametric
surfaces in the space form geometries.  Using the
method of moving frames, the proof of this classification is an
elementary exercise. 
The same proof
classifies the cyclides of Dupin in M\"obius geometry and finally in
Lie sphere geometry, where all nonumbilic Dupin immersions are Lie
sphere congruent to each other.
The idea of this proof extends to the cases of higher dimensions and greater
number of principal curvatures.  For more details and examples see the
forthcoming book by the author, Musso, and Nicolodi \cite{Surfbook}.

\section{Method of moving frames}
The method of moving frames generally refers to use of the principal
bundle of linear frames over a manifold $N$.  The text books
\cite{KobNom1, KobNom2} by
Shoshichi Kobayashi and Katsumi Nomizu, contain a seminal exposition
of this general method.  In this paper, the method of moving frames
refers to the more specialized case of the linear frames coming from a
Lie group $G$ acting transitively on a manifold $N$.  In this case, if $G_0$
denotes the isotropy subgroup of $G$ at a chosen origin $o \in N$,
then
\[
\pi:G \to N, \quad \pi(g) = g(o)
\]
is the projection map of a principal $G_0$-bundle over $N$.  A moving
frame in $N$ is any local section of this bundle.  A moving frame
along an immersed submanifold $f:M \to N$ is a smooth map $e:U \subset
M \to G$ such that $f = \pi \circ e$.  In general, there are many
local frames along $f$.  If $e:U \to G$ is a frame along $f$, then for
any smooth map $K:U \to G_0$, the map $eK:U \to G$ is also a frame
along $f$.  The method of moving frames develops a process for
reducing the local frames along $f$ to a \textit{best frame}.

The Maurer-Cartan
form of $G$ is the left-invariant $\g$-valued 1-form $\omega =
g^{-1}dg$ on $G$.  If $G$ is a matrix group, say $G$ is a closed
subgroup of the general linear group $\GL(n,\R)$, then $\omega =
(\omega^i_j) \in \gl(n,\R)$ is an $n\times n$ matrix of left invariant
1-forms on $G$ that satisfy the \textit{structure equations of $G$},
$d\omega = -\omega\wedge \omega$, which in components is
\[
d\omega^i_j = -\sum_k \omega^i_k \wedge \omega^k_j.
\]
Frame reduction is a systematic way of imposing linear relations on
the forms $\omega^i_j$ of $e^*\omega$, for a frame $e:U \to G$ along
$f$.

The following two Cartan-Darboux theorems provide the basic analytic
tools in the method of moving frames.
\begin{congruence} If $e,\tilde e:M \to G$ are smooth maps from a
  connected manifold $M$ such that $e^*\omega = \tilde e^*\omega$,
  then there exists an element $g \in G$ such that $\tilde e = ge$ on
  $M$.
\end{congruence}

\begin{existence} If $\eta$ is a $\g$-valued 1-form on a simply
  connected manifold $M$ such that $d\eta = -\eta \wedge \eta$, then
  there exists a smooth map $e:M \to G$ such that $\eta = e^*\omega$.
\end{existence}

\section{Euclidean space} Consider an immersion $\x:M \to \R^3$, where
$M$ is always now a connected surface.  Our group now is the Euclidean
group $\E(3) = \R^3 \rtimes \SO(3)$, which is represented as a matrix
group by $(\y,A) = \begin{pmatrix}1&0 \\\y & A \end{pmatrix} \in
\GL(4,\R)$.  It acts transitively by $(\y,A)\x = \y + A\x$.  The
principal $\SO(3)$-bundle $\pi:\E(3) \to \R^3$, given by $(\y,A)\0 =
\y$, is the bundle of all oriented orthonormal frames on $\R^3$.
A frame along $\x$ is a smooth map $(\x,e):U \subset M \to \E(3)$.  If
$\ee_i$ denotes column $i$ of $e \in \SO(3)$, then $e =
(\ee_1,\ee_2,\ee_3)$ is an orthonormal frame at each point of $\x$.
The pull-back of the Maurer-Cartan form $(\x,e)^{-1}d(\x,e) =
(\theta,\omega)$ satisfies
\[
d\x = \sum_1^3 \theta^i\ee_i, \quad d\ee_i = \sum_1^3 \omega^j_i
\ee_j,
\]
where $\theta = (\theta^i)$ and $\omega = (\omega^i_j) = -\t\omega$.
The structure equations are
\[
d\theta^i = -\sum_1^3\omega^i_j\wedge \theta^j, \quad d\omega^i_j = - \sum_1^3
\omega^i_k\wedge \omega^k_j.
\]
The best frame along $\x$ is one for which $\ee_3$ is
normal to $\x$ and $\ee_1$ and $\ee_2$ are principal directions.
The pull-back of the Maurer-Cartan form, $(\x,e)^{-1}d(\x,e) =
(\theta,\omega)$, expresses these conditions by
\[
\theta^3 = 0, \quad \omega^3_1 = a \theta^1, \quad \omega^3_2 =
c\theta^2,
\]
where the functions $a,c:U \to \R$ are the principal curvatures.
The immersion is totally umbilic if $a=c$ at every point of $M$, in
which case $a$ is
constant on $M$, and $\x(M)$ lies in a plane, if $a = 0$ or $\x(M)$
lies in a sphere of radius $1/|a|$, if $a \neq 0$.  

\begin{definition} The immersion $\x:M \to \R^3$ is
  \textit{isoparametric} if its principal curvatures are constant on $M$.
The immersion is \textit{Dupin} if
  its principal curvatures are distinct and each is constant along its
  own lines of curvature.
\end{definition}

Nonumbilic isoparametric immersions are Dupin.
The proof of the
following elementary theorem for nonumbilic isoparametric immersions
is the prototype of the classification of
Dupin hypersurfaces in spheres.

\begin{theorem} If $\x:M \to \R^3$ is isoparametric with distinct
  principal curvatures
$a$ and $c$, then $ac = 0$,
  and $\x(M)$ lies in a circular cylinder of radius $R$, where $1/R$
  is the absolute value of the nonzero principal curvature.
\end{theorem}

\proof  For the frame $(\x,e):U \to \E(3)$ above, 
\[
\theta^3 = 0, \quad \omega^3_1 = a\theta^1,\quad \omega^3_2 =
c\theta^2,
\]
and $\theta^1\wedge \theta^2$ is never zero, since $\x$ is an
immersion.  Then $a$ and $c$ constant and distinct combined with the
structure equations, implies $ac = 0$ and $\omega^1_2 = 0$.  Suppose
$a \neq 0$.  Now the pull back of the Maurer-Cartan form satisfies
\[
\theta^3 = 0, \quad \omega^1_2 = 0,\quad \omega^3_1 = a\theta^1,\quad
\omega^3_2 = 0.
\]
These equations on the components of the Maurer-Cartan form itself
define a completely integrable, left-invariant,
2-plane distribution $\h$ on 
$\E(3)$, so $\h$ is a Lie subalgebra of
$\mathcal E(3)$.  If $H$ is the connected Lie subgroup whose Lie
algebra is $\h$, then the maximal integral surfaces of this
distribution are the right cosets of $H$.  Any frame of the above type
along $\x$ is an integral surface of $\h$, so it lies in a right
coset $(\y,A)H$, and the projection $\x(U) = (\y,A)H\0$.  $M$
connected implies all of $\x(M)$ lies in this projection.  But this
projection is congruent to $H\0$, which by exponentiation of $\h$ is
the cylinder $x^2 + (z-\frac1a)^2 = \frac1{a^2}$ in $\R^3$.
\qed

Relative to the best frame $(\x,e):U \to \E(3)$ above, $da =
a_1\theta^1+ a_2 \theta^2$ and the lines of curvature of $a$ are the
integral curves of $\theta^2 = 0$, so $a$ is constant along its lines
of curvature iff $a_1 = 0$.  Similarly, $c_2=0$ is the condition for
$c$ to be constant along its lines of curvature.  A Dupin immersion is
a solution of the system of PDE
\[
a_2 = p(a-c), \quad c_1= q(a-c),\quad p_2 - q_1 = ac + p^2 + q^2,
\]
where $da = a_2 \theta^2$, $dc = c_1 \theta^1$, and 
$\omega^2_1 = p\theta^1 + q\theta^2$.  The situation is
essentially the same in the three space form geometries.

\section{Space form geometries} The space form geometries are
\begin{itemize}
\item Euclidean space $\R^3 = \E(3)/\SO(3)$.
\item Spherical geometry $\sph = \SO(4)/\SO(3) = \{\x=\sum_0^3 x^i\e_i
  \in \R^4:\x\cdot \x =  \sum_0^3 (x^i)^2 = 1\}$.
\item Hyperbolic geometry $\H^3 = \SO(3,1)/\SO(3) = \{\x= \sum_1^4
  x^i\e_i \in \R^{3,1}: \langle \x,\x \rangle = \sum_1^3 (x^i)^2 - (x^4)^2 = 0\}$.
\end{itemize}

Let
$
\span\{\e_0,\e_1,\e_2,\e_3,\e_4,\e_5\} = \R^{4,2}
$
for the standard orthonormal basis of 
signature $++++--$, so
\[
\aligned
\R^3 &= \span\{\e_1,\e_2,\e_3\}, \quad \R^4 = \span\{\e_0,\e_1,\e_2,\e_3\}, \\
\R^{3,1} &= \span\{\e_1,\e_2,\e_3,\e_4\}, \quad \R^{4,1} =
\span\{\e_0,\e_1,\e_2,\e_3,\e_4\}.
\endaligned
\]

\begin{theorem}
The nonumbilic isoparametric immersions in $\sph$ have constant
principal curvatures $a \neq c$ satisfying $ac +1=0$.  They
are the 1-parameter family of circular
tori $\mathbf S^1(r) \times \mathbf S^1(s) \subset \sph$, where $r =
\cos \alpha$, $s = \sin \alpha$, and the parameter is $0<\alpha
\leq \pi/4$. Here $\mathbf
S^1(r)$ is the circle in $\R^2$ with center at the origin and radius
$r$, and $a = -\tan \alpha$, $c = \cot \alpha$.  
\end{theorem}

\begin{theorem}
The nonumbilic isoparametric immersions in $\H^3$ have constant
principal curvatures $a \neq c$ satisfying $ac -1 = 0$.  They
are the 1-parameter family of
circular hyperboloids $\mathbf S^1(\frac ab)\times \H^1(\frac1b) \in
\H^3$, where the parameter is $0<a<1$ and $b = \sqrt{1-a^2}$.  Here
$\mathbf S^1(\frac ab) \subset \R^2 = \span\{\e_1,\e_2\}$, and
$\H^1(\frac1b) \subset \R^{1,1} = \span\{\e_3,\e_4\}$ is the hyperboloid
$z^2-w^2=-b^2$.
\end{theorem}

The moving frame proofs in both cases are essentially identical to
that for the Euclidean case.
Dupin immersions into $\sph$ and $\H^3$ are defined in the same way as
in the Euclidean case.

\begin{figure}
\includegraphics[scale=0.5]{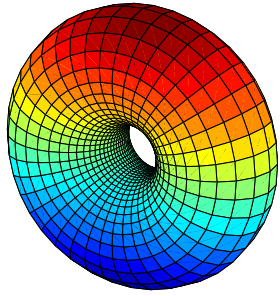}
\caption{Stereographic projection of circular torus with $\alpha = \pi/4$.}
\label{Fig1}
\end{figure}

The space form geometries are related by conformal diffeomorphisms that
send spheres to spheres.  These are, stereographic projection
\[
\mathcal S:\sph\setminus\{-\e_0\} \to \R^3, \quad \mathcal S(\sum_0^3
x^i\e_i) = \frac{\sum_1^3 x^i\e_i}{1+\x^0},
\]
with inverse
\[
\mathcal S^{-1}(\sum_1^3 y^i \e_i) = \frac{(1-\sum_1^3 (y^i)^2)\e_0 +
  2\sum_1^3y^i\e_i}{1+\sum_1^3 (y^i)^2}.
\]
Stereographic projection of the circular torus with parameter $\alpha =
\pi/4$ is shown in Figure~\ref{Fig1}. It is a Dupin immersion, but not
isoparametric, into $\R^3$.

\begin{figure}
\includegraphics[scale=0.5]{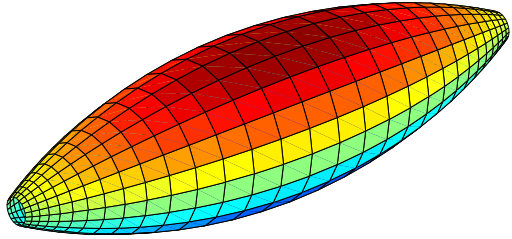} 
\caption{Hyperbolic stereographic projection of circular hyperboloid
  with $a = 1/2$.}
\label{Fig2}
\end{figure}

Hyperbolic stereographic projection onto the unit ball
$\B^3\subset \R^3$ is
\[
\mathfrak s:\H^3 \to \B^3, \quad \mathfrak s(\sum_1^4 x^i\e_i) =
\frac{\sum_1^3 x^i\e_i}{1+\x^4},
\]
with inverse
\[
\mathfrak s^{-1}(\sum_1^3 y^i\e_i) = \frac{2\sum_1^3y^i\e_i+
  (1+\sum_1^3 (y^i)^2)\e_4}{1-\sum_1^3 (y^i)^2}.
\]
It is an isometry onto the ball with the Poincar\'e metric $I_\B =
\frac{4d\y\cdot d\y}{(1-|\y|^2)^2}$, and it is a conformal embedding
when regarded as a map into $\R^3$ by the conformal inclusion $\B^3
\subset \R^3$. 
Then $\mathcal S^{-1} \circ \mathfrak s:\H^3 \to \sph$ is a conformal
embedding.  The projection by $\mathfrak s$ of the circular
hyperboloid with parameter 
$a = 1/2$ is shown in Figure~\ref{Fig2}.  It is an isoparametric
immersion into the Poincar\'e ball.  It is Dupin, but 
not isoparametric, into $\R^3$.

\section{Tangent spheres}  These conformal diffeomorphisms do not send
isoparametric immersions 
to isoparametric immersions, but they do send Dupin immersions to
Dupin immersions.  This is easily seen if we reformulate the Dupin
condition in terms of tangent sphere maps.  The reformulated
definition makes sense in M\"obius and Lie sphere geometries, as well
as in the space forms.

\begin{definition} The \textit{oriented sphere of center $\m \in \sph$
    and signed radius $r$}, where $(0<r<\pi)$, is
\[
S_r(\m) = \{\x \in \sph: \x\cdot \m = \cos r\}
\]
oriented by the unit normal $\n(\x) = \frac{\m -\cos r\ \x}{\sin r}$.
Note that $S_r(-\m)$ is $S_{\pi-r}(\m)$ with the opposite orientation.
\end{definition}

The set of all oriented spheres in $\sph$ is identified with the
smooth hypersurface 
\[
\Os = \{ S = \sum_0^4 s^i\e_i \in \R^{4,1}: \langle S,S \rangle =
\sum_0^3 (s^i)^2 - (s^4)^2 = 1\},
\]
by $S_r(\m) \leftrightarrow \frac{\m + \cos r\ \e_4}{\sin r}$ and $S =
\sum_0^4 s^i\e_i \leftrightarrow S_r(\frac{\sum_0^3 s^i\e_i}{\sin
  r})$, where $\cot r = s^4$.

Let $\x:M \to \sph$ be an immersed surface with
unit normal vector field $\ee_3$ along $\x$.

\begin{definition} An \textit{oriented tangent sphere} to $\x$ at a
  point is an oriented 
  sphere tangent to the surface at $\x$ with unit normal at $\x$ equal
  to $\ee_3$.  It is one of the
  spheres $S_r(\cos r\ \x + \sin r\ \ee_3)$, where $0<r<\pi$.  A
  tangent sphere is a \textit{curvature sphere} to $\x$ if $\cot r$ is
  a principal curvatue of $\x$ at the point.
\end{definition}

\begin{definition} A \textit{tangent sphere map along} $\x:M \to \sph$
  is a smooth map $S:M \to \Os$ such that $S(m)$ is an oriented
  tangent sphere to $\x$ at $\x(m)$, for all $m \in M$.  It must be
  given by 
\[
S = \frac{\cos r\ (\x + \e_4) + \sin r\ \ee_3}{\sin r} = \cot r\
(\x+\e_4) + \ee_3.
\]
for smooth function $r:M \to (0,\pi)$.
\end{definition}

If $S:M \to \Os$ is a tangent sphere map along $\x$, then $dS =$
\[
((\cot r - a)\ee_1 + (\cot r)_1(\x+\e_4))\omega^1_0 + ((\cot r -
c)\ee_2 + (\cot r)_2(\x+\e_4))\omega^2_0,
\]
where $d\cot r = (\cot r)_1 \omega^1_0 + (\cot r)_2 \omega^2_0$.
Thus, $dS_m \mod (\x+\e_4)$
has rank less than 2 if and only if $\cot r$ is
a principal curvature of $\x$ at $m$ if and only if $S(m)$ is a curvature
sphere of $\x$.  If $S$ is a curvature sphere map along $\x$, then
$dS$ has rank less than 2 at every point if and only if $\cot r$ is
constant along its lines of curvature.  This gives us an alternate,
but equivalent, definition of Dupin immersion into $\sph$.

\begin{definition} An immersion $\x:M \to \sph$ is
Dupin if $dS$ is singular at every point of $M$,
for any curvature sphere map $S:M \to \Os$ along $\x$.
\end{definition}

\section{M\"obius geometry}  M\"obius space is $\sph$ acted upon
by the group of all of its conformal diffeomorphisms. With the projective
description
\[
\M = \{[q] = [\sum_0^4q^i\e_i]\in P(\R^{4,1}):\sum_0^3(q^i)^2 -
(q^4)^2=\langle q,q \rangle = 0\},
\]
we have a conformal diffeomorphism
\[
f_+:\sph \to \M, \quad f_+(\sum_0^3x^i\e_i) = [\sum_0^3 x^i\e_i +\e_4],
\]
and the group of all conformal transformations on $\sph$ is
represented on $\M$ as the group of linear transformations
$\SO(\R^{4,1})$.  Let $\Mob \subset \GL(5,\R)$ be the matrix 
representation of this group in the
basis 
\[
\d_0 = \frac{\e_4+\e_0}{\sqrt2},\quad \d_i = \e_i,\quad \d_4 =
\frac{\e_4 - \e_0}{\sqrt2}
\]
of $\R^{4,1}$.  Then $\M = \Mob/G_0$,
where $G_0$ is the isotropy subgroup at the origin $[\d_0] =
f_+(\e_0)$.  The columns $\Y_a$ of an element $Y \in \Mob$ form a
M\"obius frame of $\R^{4,1}$:
\[
(\langle \Y_a,\Y_b \rangle)_{0\leq a,b \leq 4} = g = \begin{pmatrix}
0&0& -1 \\ 0&I_3 &0 \\ -1&0&0 \end{pmatrix}.
\]
The space forms are conformally contained in
$\M$ by
\[
f_+:\sph \to \M, \quad f_0 = f_+\circ \mathcal S^{-1}:\R^3 \to \M,
\quad f_- = f_0\circ \mathfrak s:\H^3 \to \M.
\]
These embeddings are equivariant with natural monomorphisms 
\[
F_+:\SO(4) \to \Mob, \quad F_0:\E(3) \to \Mob, \quad F_-:\SO(3,1) \to
\Mob.
\]
A M\"obius frame along an immersed surface $f:M \to \M$ is a
smooth map $Y:U \subset M \to \Mob$ such that $f = [\Y_0]$ on $U$.  
Then $d\Y_b = \sum_0^4 \omega^a_b \Y_a$ and
$\omega^4_0 = -\langle d\Y_0, \Y_0\rangle = 0$, since $\langle
\Y_0,\Y_0 \rangle = 0$ is constant on $U$.
It
is a first order frame if $\omega^3_0 = 0$ on $U$.  That is, $d\Y_0 =
\sum_0^2 \omega^a_0 \Y_a$.
For such a frame, $\Y_3:U
\to \Os$ is a tangent sphere map along $f$, since $\langle \Y_0,\Y_3
\rangle = 0$ implies that $[\Y_0]$ lies on this sphere, and $\langle
d\Y_0, \Y_3\rangle =0$ implies that the sphere is tangent to $f$ at
each point.  Moreover, for any smooth function $r:U \to \R$, the map
$\Y_3 + r\Y_0:U \to \Os$ is a tangent sphere map along $f$.  This is a
curvature sphere map iff
\[
d(\Y_3 + r\Y_0) \equiv \sum_1^2 (\omega^a_3+ r \omega^a_0)\Y_a \mod \Y_0 
\]
has rank less than 2 at each point of $U$, which is equivalent to
\[
(\omega^1_3 + r \omega^1_0)\wedge (\omega^2_3 + r\omega^2_0) = 0
\]
at every point of $U$.  If we assume distinct curvature spheres at
each point of $M$, then about each point we can reduce to a best
frame $Y:U \to \Mob$, characterized by
\[
\aligned
\o^3_0 &= 0, \; \omega^1_0\wedge \omega^2_0 >0 &\quad &\mbox{ (first order)} \\
\o^3_1 - i \o^3_2 &= \o^1_0 + i \o^2_0 &\quad &\mbox{ (second order)} \\
\o^0_3 &=0 &\quad &\mbox{ (third order).}
\endaligned
\]
The structure equations then imply
\[
\aligned
\omega^2_1 &= q_1 \omega^1_0 + q_2 \omega^2_0, \quad
\omega^0_0 = -2(q_2 \omega^1_0 - q_1 \omega^2_0), \\
\omega^0_1 &= p_1 \omega^1_0 + p_2 \omega^2_0, \quad \omega^0_2 = -p_2
\omega^1_0 + p_3 \omega^2_0,
\endaligned
\]
for smooth functions $q_1, q_2, p_1, p_2, p_3:U \to \R$, which satisfy
\[
\aligned
d(q_2 + iq_1)\wedge \varphi &= -\frac12
(p_1 + p_3 +1+q_1^2 + q_2^2 +ip_2)\varphi \wedge \bar\varphi, \\ 
d(p_1 + p_3 -i 2p_2) &\wedge \varphi + d(p_1-p_3) \wedge \bar\varphi = 
(2p_2 + i(p_1+p_3))(q_1+iq_2)\varphi \wedge \varphi,
\endaligned
\]
where $\varphi = \omega^1_0 + i \omega^2_0$.  Relative to a best
frame $Y:U \to \Mob$, the curvature sphere maps are $\Y_3 +\epsilon
\Y_0$, where $\epsilon = \pm1$.  The Dupin condition is that these
maps are singular at every point.  But
\[
d(\Y_3 + \epsilon \Y_0) = (\epsilon -1)\omega^1_0 \Y_1 + (\epsilon
+1)\omega^2_0 \Y_2 + \epsilon \omega^0_0 \Y_0
\]
is singular at each point, for both choices of $\epsilon = \pm1$ iff
$\omega^0_0 = 0$ iff $q_1 = q_2 = 0$ on $U$.  From the structure
equations it then follows that $f:M \to \M$ is Dupin iff $p_2 = 0$ as
well,
\[
\omega^0_0 = 0 = \omega^2_1, \quad \omega^0_1 = p_1 \omega^1_0, \quad
\omega^0_2 = p_3 \omega^2_0,
\]
and
\[
p_1 + p_3 = -1, \quad p_1-p_3= 2C,
\]
for some constant $C \in \R$.  Hence, a best frame field 
$Y:U \to \Mob$ along a Dupin immersion is an integral
surface of the completely integrable, left-invariant 2-plane
distribution on $\Mob$ defined, for each $C \in \R$, by the Lie
subalgebra $\h_C \subset \mob$ 
defined by the equations
\[
\aligned
\omega^0_0 &= 0, \quad \omega^0_1 = (-\frac12 +C)\omega^1_0, \quad \omega^0_2
= (-\frac12 -C)\omega^2_0, \quad \omega^0_3 = 0, \\
\omega^2_1 &= 0,\quad \omega^3_1 = \omega^1_0, \quad \omega^3_2 =
-\omega^2_0,\quad \omega^3_0 = 0.
\endaligned
\]
If $H_C\subset \Mob$ is the connected Lie subgroup of $\h_C$, then the
maximal integral submanifolds of this distribution are the right
cosets of $H_C$.  The projection $H_C[\d_0]$ is a Dupin surface in
$\M$, with $\pm C$ essentially the same.  These are:
\begin{itemize}
\item $f_+$ of circular tori of parameter $0<\alpha\leq \pi/4$ for
  $0\leq C= \cos 2\alpha <1$.
\item $f_0$ of the circular cylinder of radius 1 for $C = 1$.
\item $f_-$ of the circular hyperboloids of parameter $0<a<1$ for $C=
  \frac{a^{-1}+a}{a^{-1}-a} >1$.
\end{itemize}
We have thus proved that any Dupin immersion of a connected surface
into Euclidean space is contained in a conformal transformation of a
nonumbilic isoparametric immersion into some space form.

\section{Lie sphere geometry}  The fundamental role of tangent sphere
maps along an immersion $f:M \to \M$ provides motivation to
linearize the set of all oriented tangent spheres at a point of $f$.
We call such a set
a \textit{pencil of oriented spheres} at a point of $f$.  Projectivize
$\Os$ by the \textit{Lie quadric}
\[
Q = \{[q] \in P(\R^{4,2}): \langle q,q \rangle = 0\},
\]
by the inclusion $\Os \subset Q$ given by $S \mapsto S + \e_5$.
M\"obius space itself is a subset of $Q$ by the natural inclusion
$\R^{4,1} \subset \R^{4,2}$.  Then
\[
Q = \M \cup \Os = \{[q]:\langle q,\e_5 \rangle = 0\} \cup
\{[q]:\langle q,\e_5 \rangle \neq 0\}
\]
is a disjoint union comprising all oriented spheres in $\M$, including
the point spheres $\M$ (radius is zero), which have no orientation.  A
pencil of oriented spheres is a line in $Q$, which is the linear span
of any two orthogonal elements $[S_0], [S_1] \in Q$, denoted $[S_0,S_1]$.  
\begin{definition} The space of Lie sphere geometry is $\Lambda$, the set of
all lines in $Q$.  
\end{definition}
$\Lambda$ is a smooth manifold of dimension five.  Any
point $\lambda \in \Lambda$ contains a unique point sphere, which then
defines the \textit{spherical projection map} 
\[
\sigma:\Lambda \to \M, \quad \sigma(\lambda) = \lambda \cap \M,
\]
which is a fiber 
bundle with standard fiber equal to $\s2$.  The group $\SO(\R^{4,2})$
preserves $Q$ and acts transitively on $\Lambda$.  As origin of
$\Lambda$ choose
\[
o= [\l_0, \l_1], \quad \l_0 = \frac{\e_5 + \e_0}{\sqrt 2}, \quad \l_1 =
\frac{\e_4+\e_1}{\sqrt2} \in \R^{4,2}.
\]
A \textit{Lie frame} of $\R^{4,2}$ is a basis $\Tb_0,\dots,\Tb_5$ for
  which
\[
(\langle \Tb_a,\Tb_b \rangle)_{0\leq a,b \leq 5} = \hat g =
\begin{pmatrix} 0&0& -L \\
0& I_2 & 0 \\ -L & 0&0 \end{pmatrix},
\]
where $L = \begin{pmatrix}0&1 \\ 1&0 \end{pmatrix}$.  The \textit{Lie
  sphere group} is 
\[
G =\{ T \in \GL(6,\R): \t T \hat g T = \hat g \},
\]
the
representation of $\SO(\R^{4,2})$ in the Lie frame $\l_0,\dots,\l_5$, where
\[
\l_2 = \e_2, \quad \l_3 = \e_3,\quad \l_4 =
\frac{\e_4-\e_1}{\sqrt2},\quad \l_5 = \frac{\e_5-\e_0}{\sqrt2}.
\]
Then $T \in G$ iff the columns $\Tb_0,\dots,\Tb_5$ of $T$ form a Lie
frame.  Its Lie algebra in this representation is
\[
\g = \{\mathcal T \in \gl(6,\R): \t \mathcal T \hat g + \hat g
\mathcal T = 0\}.
\]
The Maurer-Cartan form of $G$ is $\omega = (\omega^a_b) \in \g$.
If $G_0$ is the isotropy subgroup of $G$ at $o$, then the
projection
\[
\pi:G \to \Lambda, \quad \pi(T) = [\Tb_0,\Tb_1]
\]
is a principal $G_0$-bundle.  A contact structure is given on
$\Lambda$ by the 1-forms $\omega^4_0 = -\langle d\Tb_0,\Tb_1 \rangle$
for any local section $T:U\subset \Lambda \to G$.  Lie sphere geometry
is the study of Legendre immersions $\lambda:M \to \Lambda$.  Any
immersion $f:M \to \M$ with tangent sphere map $S:M \to \Os$ along it
has a \textit{Legendre lift} 
\begin{equation}\label{LL}
\lambda = [f, S+\e_5]:M \to \Lambda.
\end{equation}
This is a Legendre immersion, for if $F:U \subset M\to \R^{4,1}$ is
any lift of $f$, then $\langle F,S+\e_5 \rangle = \langle F,S \rangle
= 0$, so $[F,S+\e_5] \in \Lambda$, and
the pull-back by $\lambda$ of the contact
structure is $-\langle dF,S+\e_5 \rangle = -\langle dF,S\rangle = 0$,
all since $S$ is a tangent sphere map along $f$.  The spherical
projection of the Legendre lift of $f$ is again $f$.
The spherical projection of a general Legendre immersion $\lambda$ is
a smooth map $\pi\circ\lambda:M \to \M$, but not necessarily an
immersion.  

\begin{example}\label{example}
$\lambda:\R^2 \to \Lambda$, $\lambda(u,v) = [\sb_0(u),\sb_1(v)]$,
where
\begin{equation*}\label{eq:Li:89.1b}
\sb_0(u) = \cos u \e_0 + \sin u \e_3 + \e_4, \quad
\sb_1(v) = \cos v \e_1 + \sin v \e_2 + \e_5.
\end{equation*}
is a smooth Legendre immersion.  Its spherical projection $\sigma\circ
\lambda:\R^2 \to \M$ is $\sigma \circ \lambda(u,v) = [\sb_0(u)]$,
which is $f_+$ of the great circle $\cos u\ \e_0 + \sin u\ \e_3$ in
$\sph$.  In particular, its spherical projection is singular at every
point of $M = \R^2$.
\end{example}

Given a Legendre immersion $\lambda:M \to \Lambda$,
for each $m\in M$, the line $\lambda(m)$ is the pencil of tangent
spheres at $m$.
A smooth map $S:U \subset M \to Q$ such that
$S(m) \in \lambda(m)$, for each $m\in U$ is a \textit{tangent sphere
  map} along $\lambda$.  The definitions of curvature sphere map and
Dupin Legendre immersions is the same as for immersions into $\M$.
The important, but elementary, fact is that the Legendre lift~\eqref{LL}
of a
Dupin immersion $f:M \to \M$ with tangent sphere map $S:M \to \Os$
is a Dupin Legendre immersion.

An elementary calculation verifies that the Legendre immersion of
Example~\ref{example} is Dupin.

The best Lie frame field $T:U \to G$ along a
Dupin Legendre immersion $\lambda:M \to \Lambda$ is characterized by
$[\Tb_0]=S_0$ and $[\Tb_1]=S_1$ are the curvature spheres of $\lambda$ and
\begin{equation*}
\aligned
&\text{Order 1: }\; \omega^2_0 = 0 = \omega^3_1, \; \theta^2 = \omega^2_1,\;
\theta^3 = \omega^3_0,\; \theta^2 \wedge \theta^3 \neq 0, \\
&\text{Order 2: }\; 0=\omega^1_0 = \omega^0_1 =\omega^2_3  = -\omega^3_2,
\\
&\text{Order 3: }\; 0=\omega^0_2  = \omega^1_3=\omega^0_4.
\endaligned
\end{equation*}
By the structure equations of $G$,
\[
d\theta^2 = p\theta^2\wedge \theta^3, \quad d\theta^3 = q
\theta^2\wedge \theta^3,
\]
for smooth functions $p,q:U \to \R$.  The remaining entries of
$\omega$ are given by
\begin{equation*}
\omega^0_0 = q\theta^2 + t \theta^3, \quad \omega^1_1 = u \theta^2 - p
\theta^3,\quad \omega^0_3 = c_i\theta^i, \quad \omega^1_2 =
d_i\theta^i,
\end{equation*}
for smooth functions $t,u, c_i, d_i:U \to \R$, where $i = 2,3$.
Taking the exterior differential of these forms, we get
\begin{equation*}
\aligned
dq \wedge \theta^2 + dt \wedge \theta^3 &=-(c_2 + q(p+t))\theta^2\wedge
\theta^3, \\ 
du \wedge \theta^2 - dp \wedge \theta^3 &= (d_3 +
p(q-u))\theta^2\wedge \theta^3.
\endaligned
\end{equation*}

\begin{figure}
\includegraphics[scale=0.5]{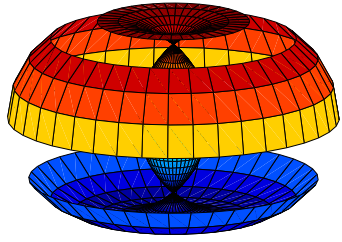}
\caption{Spherical projection followed by hyperbolic stereographic
  projection of a right coset of $H$.}
\label{Fig7}
\end{figure}

\begin{lemma}\label{le:Li:70}  
The left-invariant
6-plane distribution $\mathcal D^\perp$ defined on $G$ by
\[
\mathcal D = \{\omega^2_0, \omega^3_1, \omega^1_0, \omega^0_1,
\omega^2_3, \omega^0_2, \omega^1_3, \omega^0_4, \omega^4_0\}
\]
is completely integrable.  Its maximal integrable manifolds are the right
cosets of the connected 6-dimensional Lie subgroup $H$ of $G$ whose Lie algebra
is
\begin{equation*}\label{eq:Li:92}
\mathfrak h = \mathcal D^\perp = \{ \mathcal T \in \g: \alpha(\mathcal
T) = 0, \text{ for all $\alpha \in \mathcal D$.}\}.
\end{equation*}
\end{lemma}

\begin{theorem} If $\lambda:M \to \Lambda$ is any Dupin Legendre
  immersion, then $\lambda(M)$ is
  an open submanifold of $AHo$, for some element $A \in G$.  
\end{theorem}

\proof Each point of $M$ has
  a neighborhood $U$ on which there exists a best Lie frame field $T:U \to G$
  along $\lambda$ with $T^{-1}dT \in \h$ on $U$, where $\h \subset \g$
  is the Lie algebra of $H$.  Thus, $T:U \to G$ is an integral surface
  of the 6-plane distribution defined on $G$ by $\h$, so $T(U) \subset
  AH$, for some element $A \in G$, since the integral submanifolds of
  $\h$ are the right cosets of $H$.  If $M$ is connected, then $S(V)
  \subset AH$, for any best Lie frame $S:V\subset M \to G$.  Hence,
  $\lambda(M) \subset AHo$.
\qed

Calculation of the best Lie sphere frame field along the Dupin
Legendre immersion in Example~\ref{example}
reveals that $\lambda(\R^2) = Ho$.  Thus, the spherical projection
$\sigma(Ho)$ is a great circle.  For fixed $t\in \R$, if $A \in G$
is the matrix of 
\[
\begin{pmatrix} \cosh t & 0 & \sinh t \\
0 & I_4 & 0 \\
\sinh t & 0 & \cosh t \end{pmatrix} \in \SO(4,2),
\]
in the Lie frame, then $\mathcal S \circ f_+^{-1}\circ \sigma (AHo)$,
the stereographic projection of the spherical projection of $AHo$,
is a surface of revolution obtained by rotating a circle about an
axis that intersects it.  It is singular at the points of intersection.
Figure~\ref{Fig7} shows it with part of the
outer surface removed to reveal the inner part and the singularities.

In conclusion, every connected Dupin immersion in $\R^3$ is obtained
in this way for some choice of element $A \in G$.

\bibliography{Bibliography}
\bibliographystyle{alpha}

\end{document}